\documentstyle[12pt]{article}
\title{EXPECTED NUMBER OF LOCAL MAXIMA OF  SOME GAUSSIAN   RANDOM  POLYNOMIALS}
\author{S. Shemehsavar, S. Rezakhah\thanks{Faculty of
 Mathematics and Computer Sciences, Amirkabir
 University of Technology,Tehran, Iran, email: shemehsavar@aut.ac.ir, email: rezakhah@aut.ac.ir} }
\date{}

\begin{document}

\maketitle

\pagestyle{plain}

\renewcommand{\theequation}{\arabic{section}.\arabic{equation}}

\begin{abstract}
 Let $Q_n(x)=\sum _{i=0}^{n} A_{i}x^{i}$ be a random algebraic polynomial where
 the coefficients
 $A_0,A_1,\cdots $ form a sequence of centered Gaussian random variables. Moreover,
  assume that
the increments $\Delta_j=A_j-A_{j-1}$, $j=0,1,2,\cdots$ are
independent, $A_{-1}=0$. The coefficients can be considered as $n$
consecutive observations of a Brownian motion. We study  the
asymptotic behaviour of the expected number of local maxima of
$Q_n(x)$ below level $u=O(n^k)$, for some  $k>0$.

\end{abstract}

AMS(2000) subject classifications. Primary 60H42, Secondary 60G99.

 \setcounter{equation}{0}
\newcommand{\be}{\begin{equation}} \newcommand{\ee}{\end{equation}}
\section{Introduction}
\hspace{.2in}
 The theory of the expected number of real zeros of
random algebraic polynomials was addressed in  the fundamental
work of M. Kac\cite{kac1} (1943).
 The works  Wilkins \cite{wil}, and Farahmand
 \cite{far1}, \cite{far3} and Sambandham \cite{sam1,sam2} are other fundamental
 contributions  to the subject. For various aspects on random
 polynomials see Bharucha-Reid and Sambandham \cite{bhs}, and
 Farahmand\cite{far2}.

 There has been recent interest in cases where the
coefficients form certain random processes, Rezakhah and
Shemehsavar \cite{rs3},
  Rezakhah and Soltani \cite{rs1,rs2}.

Let  $A_{0},A_{1},\cdots  $ be a mean zero Gaussian random
sequence for which the increments  $\Delta _i =A_i -A_{i-1}, \;
i=1,2,\cdots $ are independent, $A_{-1}=0$. The sequence  $A_0,A_1
\cdots $
 may be considered as successive Brownian points,
 i.e., $A_j=W(t_j),\; j=0,1,\cdots $, where $t_0<t_1< \cdots $
and $\{W(t),\;t\geq 0\}$ is the standard Brownian motion.
In this physical interpretation, Var$(\Delta_j)$ is the distance
 between successive times $t_{j-1},\; t_j$.
  Let
$$ Q_{n}(x)=\sum_{i=0}^{n}A_{i}x^{i},\, \; -\infty<x<\infty , \eqno{(1.1)}$$
  We note that $A_j=\Delta_0 +\Delta _1
+\cdots +\Delta_j, \;\;\; j=0,1,\cdots $, where $ \Delta_i \sim
N(0,\sigma_i^2 )$ and $\Delta_i$ are independent, $i=0,1, \cdots
$.
  Thus $
Q_n(x)=\sum_{k=0}^n (\sum _{j=k}^n x^j)\Delta_k=\sum_{k=0}^n
a_k(x)\Delta _k, $, $\;
  Q'_n(x)=\sum_{k=0}^n b_k (x)\Delta _k,
$, and $Q_n''(x)=\sum_{k=0}^n d_k (x)\Delta _k,$ where
    $$a_k (x) = \sum _{j=k}^n x^j,\;\;
 \;  b_k (x) = \sum _{j=k}^n
jx^{j-1},\;\; d_k (x) = \sum _{j=k}^n j(j-1)x^{j-2} \hspace{.2in}
k=0,\cdots , n.  \eqno{(1.2)} $$

In this paper we study the asymptotic behavior of the expected
number of  local maximas of $Q_n(x)$.
  We say $Q_n(x)$ has a local maxima at $t=t_i$ if $Q'_n(x)$ has a
  down-crossing of the level zero at $t_i$. A local maxima which
  we consider here, is a maxima  that occurs when
  $Q_n(x)$ is   below  level $u$.   The total number of
  down-crossing of the level zero by $Q'_n(x)$ in $(a,b)$ is
  defined as $M(a,b)$, and these occur at the points
  $a<t_1<t_2<\cdots <t_{M(a,b)}<b$. We define $M_u(a,b)$ as the
  number of  zero-down crossing    by $Q'_n(x)$ at those points $t_i\in (a,\,b)$, where   $Q(t_i)\leq u$.

Rice [1945, pp 71] showed that for any function of the random
variables $A_0, A_,\cdots A_n$ and $x$, say here $U=Q_n(x)$, the
expected number of maxima of $U$ in the interval $(a,b)$ is equal
to
$$
\int_a^b
\int_{-\infty}^{\infty}\int_{-\infty}^0|t|p_x(r,0,t)dtdrdx
\eqno(1.3)$$ where  $p_x(r,s,t)$ is the joint probability density
function of   $U=Q_n(x)$, $V=\partial Q_n(x)/\partial x$, and
$W=\partial^2 Q_n(x)/\partial x^2$. Using this formula we find
that  the expected number of local maxima of $Q_n(x)$ below level
$u$, and  inside any interval $(a,\,b)$,$\,EM_u(a,b)$ is equal to
$$
EM_u(a,b)=\int_a^b
\int_{-\infty}^{u}\int_{-\infty}^0|t|p_x(r,0,t)dtdrdx \eqno(1.4)
$$
where
$$
p_x(r,0,t)=\frac{\exp(-Lr^2-2Mrt-Kt^2)}{(2\pi)^{3/2}\det(\Sigma)^{1/2}}
$$
in which $\Sigma$ is the covariance matrix of $(U,V,W)$, and
$$
K=\frac{A^2B^2-C^2}{2\det(\Sigma)},\hspace{.3in}
L=\frac{B^2D^2-F^2}{2\det(\Sigma)} \eqno(1.5)
$$
$$
M=\frac{CF-B^2E}{2\det(\Sigma)},\hspace{.3in}S=K-\frac{M^2}{4L}
$$
and
$$
\det(\Sigma)=A^2B^2D^2-A^2F^2-B^2E^2-C^2D^2+2CEF
$$
$$
A^2=\mbox{Var}(Q_n(x))=\sum_{k=1}^na_k^2(x)\sigma_k^2,\hspace{.3in}
B^2=\mbox{Var}(Q_n'(x))=\sum_{k=1}^nb_k^2(x)\sigma_k^2,$$
$$D^2\!\!=\!\mbox{Var}(Q_n''(x))\!=\!\!\sum_{k=1}^nd_k^2(x)\sigma_k^2,\hspace{.3in}
C\!\!=\!\mbox{Cov}(Q_n(x),Q_n'(x))\!=\!\!\sum_{k=1}^na_k(x)b_k(x)\sigma_k^2,$$
$$
E=\mbox{Cov}(Q_n(x),Q_n''(x))=\sum_{k=1}^na_k(x)d_k(x)\sigma_k^2,\;$$
$$
\;F=\mbox{Cov}(Q_n'(x),Q_n''(x))=\sum_{k=1}^nb_k(x)d_k(x)\sigma_k^2,
\;$$ where $ a_k(x),\; b_k(x)$ and  $d_k(x)$ are defined in (2.1).

    Using  (1.5),and the function
 $\; \mbox{erf} (t)=2\Phi(t\sqrt{2})-1,\;$,   we find
that
$$
EM_u(a,b)=\int_a^b f_n(x)dx =J_1+J_2 \eqno(1.6)
$$
$$
J_1=\frac{1}{4\pi}\int_a^bG_1\big[\mbox{erf}(G_2)+1\big]dx,\hspace{.4in}
J_2=-\frac{1}{4\pi}\int_a^b
G_1G_3\bigg[\mbox{erf}(G_4)+1\bigg]\exp(G_5)$$ where $G_1=
\bigg(2S\sqrt{2L\det(\Sigma)}\bigg)^{-1},\hspace{.1in}G_2=
u\sqrt{L},\hspace{.1in} G_4=u\sqrt{K^{-1}M^2},$
$$
\hspace{.3in}G_3=\sqrt{\frac{M^2}{LK}},\hspace{1in}G_5=-\frac{LSu^2}{K}
$$
Farahmand\cite{far3} obtained a similar formula for the case where
the coefficients  are independent, normally distributed
  with mean zero and variance one.\vspace{.04in}

\section{Asymptotic behaviour of $EM_u$ }
 In this section
we obtain the asymptotic behaviour of the expected number of local
maxima  of  $Q_{n}(x)=0$  given by (1.1).
 We prove  the following theorem for the case that
 the increments
 $\Delta_1 \cdots \Delta_n$ are independent and  have the same distribution.
 Also we assume that
  $\sigma_k^2=1$, for $k=1\cdots n$.\vspace{.03in}

  \noindent
 \textbf{Theoream(2,1)}:   Let $Q_n(x)$ be the random algebraic
 polynomial given by (1.1) for which
 $A_{j}=\Delta_{1}+...+\Delta_{j}$ where $\Delta_{i}$,$i=1,...,n$
 are independent and $\Delta_{j}\sim N(0,1)$ then the
 expected number of local maxima  of $Q_n(x)$ below level $u$  satisfies:
\begin{eqnarray*}
EM_{u}(1,\infty)=\;\;\;\hspace{.7in}\!&\!\!&\!\hspace{-.9in}\frac{0.0013074}{4\pi}+\frac{(0.0350655)u}{2(n\pi)^{3/2}}+O(n^{-1/2})
\hspace{.3in} \mbox{for}  \; u=O(n^{5/4})\\ \hspace{-.3in}
EM_{u}(0,1)\!=\,\hspace{.92in}&&\hspace{-1.22in}
\frac{2(\sqrt{35}\!-\!5)}{345\pi}\mbox{ln}\!
\bigg(\frac{n^{3/2}}{u} \bigg) \!-\! 0.001648\! -\!
\frac{(2.033388)u}{2(n\pi)^{3/2}}
   \!+\!O(n^{\!-\!1/2})
  \hspace{.1in} \mbox{for}\;  u=O(n^{5/4})\\
EM_{u}(-\infty,-1)=\hspace{.57in}&\!\!\!\!\!\!&\!\!\!\hspace{-.7in}
\frac{0.0162552}{4\pi}+\frac{(0.0997677)u}{2\pi\sqrt{n\pi}}+O(n^{-1/2})
\hspace{.3in} \mbox{for}\;
u=O(n^{1/4})\\
EM_{u}(\!-\!1,0)\!\!=\,\;\;\;\hspace{.77in}\!&\! \!&
\hspace{-1.17in} \frac{2(\sqrt{3}\!-\!1)}{11\pi}\mbox{ln}\!
\bigg(\!\, \frac{n^{1/2}}{u}\! \,\bigg) \! +\! 0.081413 \!-\!
\frac{(0.594923)u}{2\pi\sqrt{n\pi}} +O(n^{\!-\!1/2})
 \hspace{.1in}
\mbox{for} \;
 u=O(n^{\!1/4})
\end{eqnarray*}

\noindent
 \textbf{proof:} The
asymptotic behaviour is treated separately on
  the intervals $1<x<\infty$,
$-\infty<x<-1$, $0<x<1$ and $-1<x<0$.

 For  $u=O(n^{5/4})$ and $1<x<\infty$, by the
change of variable $x=1+\frac{t}{n}$ and
 the equality  $\left( 1+\frac{t}{n}
\right)^n=e^t\left(1-\frac{t^2}{n} \right) +O\left(\frac{1}{n^2}
\right).  $  Using (1.6), we find  that $$EM_u(1 ,\infty
)=\frac{1}{n} \int_{0}^{\infty}f_n\big(1+\frac{t}{n}\big)dt,\;\;\;
$$ where  by (1.5) and (1.6), and by tedious manipulation we have that
$$
n^{-1}G_{1}\left(1+\frac{t}{n}\right)=
H_{11}(t)+O(n^{-1}),\;\;G_{3}\left(1+\frac{t}{n}\right)=
H_{13}(t)+O(n^{-1}) \eqno(2.1)
$$
$$G_{2}\left(1+\frac{t}{n}\right)=\frac{2u}{n^{3/2}\sqrt{\pi}}H_{12}(t)+O(n^{-5/4}),$$
$$
 G_{4}\left(1+\frac{t}{n}\right)=\frac{2u}{n^{3/2}\sqrt{\pi}}H_{14}(t)+O(n^{-5/4}),\;\;
 G_{5}\left(1+\frac{t}{n}\right)=1+O(n^{-1/2}), $$
where
$$
H_{11}(t)={\frac {1}{192}} \bigg(-4+ \left( 32t+16+32{t}^{2}
\right) {e^{t}}\hspace{2in}$$$$+ \left( 32{t}^{
5}+208{t}^{4}+472t+124+1040{t}^{2}+736{t}^{3} \right) {e^{2t }}
$$ $$+ \left( -288+192{t}^{4}- 768t+256{t}^{3}-1152{t}^{2}
\right) {e^{3t}}$$$$\hspace{.4in}+ \left(
20-80{t}^{2}+16t+176{t}^{4}-64 {t}^{5}-704{t}^{3} \right)
{e^{4t}}$$$$\hspace{1in}+ \left( 272+224\,t+160{t}^{2} \right)
{e^{5t}}+ \left( -140+24 t \right) {e^{6t}}
 \bigg)\times $$
 $$ \bigg(35- \bigg( 32-32t+160t^2
\bigg)e^t-(294-588t-324t
^2+600t^3+216t^4+80t^5)e^{2t}$$$$\hspace{.5in}
+(544-1632t+1056t^{2}+512{t}^{3}-192{t}^{4}
)e^{3t}$$$$\hspace{1.6in}+\bigg(1012t-253-1400t^2+736t^3
-172t^4+16t^5\bigg)e^{4t}\bigg)^{1/2}\times
 $$
$$\bigg[ {\frac {115}{192}}+ \left(
-{\frac {7}{8}}-{ \frac {35}{12}}{t}^{2}+{\frac {37}{8}}t \right)
{e^{t}}\hspace{3.2in}
$$
 $$- \left( {\frac {733}{48}}+{\frac
{485}{24}}t+{\frac {523}{48}} {t}^{2}+{\frac
{1073}{24}}{t}^{3}+\frac{5}{3}{t}^{5}+{\frac {73}{12}}{ t}^{4}
\right) {e^{2t}}\hspace{1in}$$
$$+ \left( {\frac {1043}{24 }}-{\frac
{232}{3}}{t}^{4}-{\frac {125}{8}}t+{\frac {364}{3}}{t} ^{3}-{\frac
{34}{3}}{t}^{6}+162{t}^{2}-14{t}^{5} \right) {e^{3 t}}$$
$$ + \left( 43{t}^{7}\!-\!{\frac {32777}{48}}{t}^{2}\!+\!{\frac {2161}
{12}}{t}^{5}\!+\!{\frac {31}{3}}{t}^{8}\!+\!{\frac
{507}{8}}t\!+\!{\frac { 5177}{16}}{t}^{4}\!+\!{\frac
{1239}{32}}\!+\!{\frac {1949}{12}}{t}^{6}\!+\!{ \frac
{1043}{12}}{t}^{3} \right) {e^{4t}}
$$
$$+ \left( -{\frac {6887}{24}}+{\frac {7153}{24
}}t+{\frac {12731}{12}}{t}^{2}-{\frac {5627}{6}}{t}^{3}+26{t}^
{7}-{\frac {2335}{6}}{t}^{4}-{\frac {173}{3}}{t}^{6}-{\frac
{1697} {6}}{t}^{5} \right) {e^{5t}}$$$$ + \left( {\frac
{20383}{48}}+{\frac { 68}{3}}{t}^{7}\!-\!{\frac
{4023}{16}}{t}^{2}\!-\!\frac{8}{3}{t}^{8}\!-\!{\frac {7297
}{8}}t+{\frac {29851}{24}}{t}^{3}\!-\!{\frac
{3907}{24}}{t}^{4}\!-\!{ \frac {547}{6}}{t}^{6}+{\frac
{397}{4}}{t}^{5} \right) {e^{6t}}$$
$$+ \left( -{\frac {2018}{ 3}}{t}^{2}-{\frac
{2141}{8}}+{\frac {527}{2}}{t}^{4}+6{t}^{6}+{ \frac
{6749}{8}}t-{\frac {1243}{6}}{t}^{3}-{\frac {401}{6}}{t}^{ 5}
\right) {e^{7t}}$$
$$ \hspace{.6in}+
 \left( {\frac {12155}{192}}-{\frac {6281}{24}}t+{\frac {1385}{16}}
{t}^{4}+{\frac {19097}{48}}{t}^{2}+{t}^{6}-{\frac {787}{3}}{t}^{
3}-{\frac {29}{2}}{t}^{5} \right) {e^{8t}}\bigg] ^{ -1}{t}^{-1}$$
and
$$
H_{12}(t)=\bigg[-80 \left(  ( {\frac {253}{80}}-{\frac
{46}{5}}{t}^{3} +{\frac {35}{2}}{t}^{2}-1/5{t}^{5}+{\frac
{43}{20}}{t}^{4}-{ \frac {253}{20}}t )
{e^{-2t}}\right.\hspace{1in}$$$$+ ( {\frac {147}{40}}+{ \frac
{27}{10}}{t}^{4}+{t}^{5}-{\frac {81}{20}}{t}^{2}+15/2{t}^{
3}-{\frac {147}{20}}t ) {e^{-4t}}$$$$+ ( {\frac {102}{5}} t-{\frac
{34}{5}}-{\frac {66}{5}}{t}^{2}-{\frac {32}{5}}{t}^{3}+{ \frac
{12}{5}}{t}^{4} ) {e^{-3t}} -{\frac {7}{16}}{e^{-6t
}}$$$$\hspace{2in}\left.+ ( 2{t}^{2}-2/5t+2/5) {e^{-5t}} \right)
{t}^{3}$$$$\times
 \bigg(  \left( -176\,{t}^{3}-20{t}^{2}-16{t}^{5}+44{t}^{4}+4t
+5 \right) {e^{-2t}} \hspace{1in}$$$$+ \left(
184{t}^{3}+260{t}^{2}+52{t}^{4}+ 31+118t+8{t}^{5} \right)
{e^{-4t}}$$$$\hspace{1in}+ \left( 64{t}^{3}-192\,t-
288{t}^{2}-72+48{t}^{4} \right) {e^{-3t}}$$$$\hspace{1in} + \left(
56t+68+40{ t}^{2} \right) {e^{-t}}+6t-35-{e^{-6t}}+ \left(
8\,t+8{t}^{2}+4
 \right) {e^{-5t}} \bigg) ^{-1}\bigg]^{1/2}
$$

Also
$$
H_{13}(t)=\bigg(5+
 \left( 4+36t\!-\!8{t}^{2} \right) {e^{t}}\!-\!\left( 12t-24{t}^{4}-52{t}^{3}+78\!-\!150{t}^{2} \right)
{e^{2t}}\hspace{1in}$$$$\hspace{.1in}-\!\left(
120{t}^{2}\!-\!124\!-\!24{t}^{3}+156t \right) {e^{3t}}+
  \left( 132t\!-\!58{t}^{2}+8{t}^{3}\!-\!55 \right)
{e^{4t}}\bigg)$$
$$\times
\bigg[ \bigg(
 \left( 253-736{t}^{3}+172{t}^{4}+1400{t}^{2}-1012t-16{t}^{5
} \right) {e^{4t}}\hspace{1in}$$$$\hspace{.2in}+\! \left(
294\!-\!324{t}^{2}\!-\!588t\!+\!600{t}^{3}\!+\!80{t}^{5}\!+\!216{t}^{4}
\right) {e^{2\,t}}$$$$- \left( 512{t}^{3}\!+\!1056{t
}^{2}\!+\!544\!-\!1632t\!-\!192{t}^{4} \right) {e^{3t}}
\!-\!35\!+\! \left(160{t}^{2}\!+\!32\!-\!32t \right) {e^{t}}
\bigg)\hspace{.2in}$$$$\hspace{.5in}\times \bigg( (15-4t){e
^{4t}}\!-\!(32+24t){e^{3t}}\!+\!\left(
36t\!+\!18\!+\!12{t}^{2}\!+\!8{t}^{3} \right)
{e^{2t}}\!-\!8{e^{t}}t\!-\!1 \bigg) \bigg]^{-1/2}
$$
and
$$
H_{14}(t)= \left(  \left( -78+150{t}^{2}+52{t}^{3}-12t+24\,{t}^{4}
 \right) {e^{-3t}}\right. \hspace{2in}$$$$+ \left( -120{t}^{2}-156t+124+24{t}^{3}
 \right) {e^{-2t}}$$$$+ \left( 132t+8\,{t}^{3}-58{t}^{2}-55 \right)
{e^{-t}}+5\,{e^{-5t}}+ \left. \left( 36\,t+4-8{t}^{2} \right)
{e^{-4t}}
 \right) {t}^{3/2}$$$$\times \bigg((4t-15)+(32+24t)e^{-t}-(18+36t+12t^2+8t^3)e^{-2t}
 +8te^{-3t}+e^{-4t}\bigg)^{-1/2}$$$$\bigg(-35+6t+(68+56t+40t^2)e^{-t}+(5+4t-20t^2-176t^3+44t^4-16t^5)e^{-2t}$$$$
 +(64t^3-72-192t+48t^4-288t^2)e^{-3t}+(31+118t+260t^2+184t^3+8t^5+52t^4)e^{-4t}$$$$
 +(4+8t+8t^2)e^{-5t}-e^{-6t}\bigg)^{-1/2}
$$
 As $t\rightarrow \infty $ we have that
$$H_{11}(t)\sim \frac{1}{2t^{7/2}},\;\;\;\;\;\; H_{13}(t)\sim 1,\;\;\;\;\;\;
H_{12}(t)=O(t^{7/2}e^{-t}),\;\;\;\;\;\; H_{14}(t)=O(t^{7/2}e^{-t})
$$

Thus  by (2.1) and above calculations we have that
\begin{eqnarray*}EM_u(1,\infty)\!&\!=\!&\! \frac{1}{n}\int_0^{\infty}
f_n\bigg(1+\frac{t}{n}\bigg)dt\\&& =\frac{1}{4\pi} \int_0^{\infty}
H_{11}(t)dt+\frac{u}{2n^{3/2}\pi \sqrt{\pi}}
\int_0^{\infty}H_{11}(t)H_{12}(t)dt\\&& -\frac{u}{2n^{3/2}\pi
\sqrt{\pi}}\int_0^{\infty}H_{11}(t)H_{13}(t)H_{14}(t)dt
-\frac{1}{4\pi}\int_0^{\infty} H_{11}(t)H_{13}(t)dt\bigg]
\end{eqnarray*}
where$\int_0^{\infty} H_{11}(t)dt=0.02789960660$  and
$\int_0^{\infty} H_{11}(t)H_{13}(t)dt=0.02659218098$ Also
$\int_0^{\infty} H_{11}(t)H_{12}(t)dt=0.3326450540$  and
$\int_0^{\infty} H_{11}(t)H_{13}(t)H_{14}(t)dt=0.297579554$

For  $u=O(n^{1/4})$ and $-\infty<x<-1$, let $x=-1-\frac{t}{n}$
then,
  by (1.6),
$EM_u(-\infty, -1 )=\frac{1}{n} \int_{0}^{\infty}
f_n\left(-1-\frac{t}{n}\right)dt.\,$  Using (1.5), (1.6) we have
that
$$
n^{\!-1}G_{1}\left(\!-\!1\!-\!\frac{t}{n}\right)\!=\!
H_{21}(t)\!+\!O(n^{\!-1}),\;\;\;\;\;
G_{2}\left(\!-\!1\!-\!\frac{t}{n}\right)\!=\!\frac{2u}{\sqrt{n\pi}}H_{22}(t)\!+\!O(n^{-5/4}),\eqno(2.2)$$
$\;\; G_{3}\left(-1-\frac{t}{n}\right)=H_{23}(t)+O(n^{-1}), \;\;$
and
$$
 G_{4}\left(-1-\frac{t}{n}\right)= \frac{2u}{\sqrt{n\pi}}H_{24}(t)+O(n^{-5/4}),\;\;\;\;\;\;\;\;  G_{5}\left(-1-\frac{t}{n}\right)=1 +O(n^{-1/2}), $$
where
$$
H_{21}(t)= \frac{3}{2t} \left(  \left(
-3/8-2\,{t}^{5}-9/2\,{t}^{2}-2\,{t}^{3}+3/2\,{t}^{4}-3/2\,t
\right) {e^{4\,t}}\right.$$$$\left. + \left(
{t}^{5}+9/2\,{t}^{4}+3/4\,t+7\, {t}^{3}+9/2\,{t}^{2}+3/8 \right)
{e^{2\,t}}+3/4\,{e^{6\,t}}t+1/8\,{e^{ 6\,t}}-1/8 \right)$$$$
\bigg[3+(12t-6-12t^2-56t^3-56t^4-16t^5)e^{2t}+(3-12t+24t^2+32t^3-44t^4+16t^5)e^{4t}\bigg
]^{1/2}$$$$ \bigg[  \left( 7/3{t}^{8}+{\frac {221}{12}}{t}^{
5}+1/8t+{\frac {11}{32}}+{\frac {259}{48}}{t}^{4}+7/4{t}^{3}+{
\frac {35}{3}}{t}^{7}+{\frac {257}{12}}{t}^{6}+{\frac {29}{16}}{
t}^{2} \right) {e^{4\,t}}$$$$+ \left( -{\frac
{15}{8}}{t}^{3}-{\frac {55
}{48}}{t}^{2}-1/24\,t-5/4{t}^{4}-1/3{t}^{5}-{\frac {11}{48}}
 \right) {e^{2t}}$$$$+ \left( {\frac {13}{6}}{t}^{6}-1/8t-{\frac {49
}{24}}{t}^{3}-{\frac {51}{4}}{t}^{5}-8/3{t}^{8}+8/3{t}^{7}-{ \frac
{211}{24}}{t}^{4}-{\frac {11}{48}}-3/16{t}^{2} \right) {e^{6
t}}$$$$+{\frac {11}{192}}+ \left( -5/2{t}^{5}+1/24t-{\frac {23}{48}}
{t}^{2}+{t}^{6}+{\frac {13}{6}}{t}^{3}+{\frac {11}{192}}+{\frac {
25}{16}}{t}^{4} \right) {e^{8t}} \bigg] ^{-1},
 $$
$$
H_{22}(t)=2\sqrt{t}\bigg(3-(16t^5+56t^4+56t^3+12t^2-12t+6)e^{2t}
\hspace{2in}
$$$$+(16t^{5}-44t^4+32t^{3}+24t^2-12t+3)e^{4t}
\bigg)^{1/2}$$$$\times
\bigg[(1+6t)e^{6t}+(12t^4-16t^5-16t^3-36t^2-12t-3)e^{4t}+(3+6t+36t^2+56t^3+36t^4+8t^5)e^{2t}-1\bigg]^{-1/2}
$$
and $$ H_{23}(t)={\frac {1}{8}}\bigg[ \bigg(
1+(8t^4+20t^3+14t^2+4t-2)e^{2t}+(1-4t-10t^2+8t^3)e^{4t} \bigg)
^{2}$$$$\times\bigg( \left( ( -2{t}^{3}-3{t}^{2}-t-1/2
 ) {e^{2t}}+1/4+(1/4+t){e^{4t}} \right)
 $$$$\times\left(({t}^{5}-11/4{t}^{4}+3/16+2{t}^{3}+3/2{t}^{2}-3/4t
 ) {e^{4t}}+3/16 \right.$$$$\left.+( 3/4t-3/4{t}^{2}-7/2{t}^{3}-7/2
{t}^{4}-{t}^{5}-3/8 ) {e^{2t}} \right) \bigg)^{-1}\bigg]^{1/2}
,$$$$ H_{24}(t)= 2\sqrt{t}\,
\bigg(1+(8t^4+20t^3+14t^2+4t-2)e^{2t}+(1-4t-10t^2+8t^3)e^{4t}
 \bigg)$$$$\times\bigg[
 \bigg((3+6t+36t^2+56t^3+36t^4+8t^5)e^{2t}+(1+6t)e^{6t}-1$$$$+(12t^4-16t^5-16t^3-36t^2-12t-3)e^{4t}
 \bigg) \left(1-(2+4t+12t^2+8t^3)e^{2t}+(1+4t)e^{4t} \right)\bigg]^{-1/2}
$$
 As $t\rightarrow \infty $ we have that
$$H_{21}(t)\sim \frac{1}{2t^{7/2}},\;\;\;\;\; H_{23}(t)\sim 1,\;\;\;\;\;
 H_{22}(t)=O(t^{5/2}e^{-t}), \;\;\;\;\; H_{24}(t)=O(t^{5/2}e^{-t})$$

Thus by (2.2) and above calculations we have that
\begin{eqnarray*}
EM_u(-\infty,-1)\!&\!=\!&\!
\frac{1}{n}\int_0^{\infty}f_n(-1-\frac{t}{n})
\\&& =\frac{1}{4\pi }
\int_0^{\infty} H_{21}(t)dt
+\frac{u}{2\pi\sqrt{n\pi}}\int_0^{\infty}H_{21}(t)H_{22}(t)dt
\\&&-\frac{u}{2\pi \sqrt{n\pi}}\int_0^{\infty}  H_{21}(t)H_{23}(t)H_{24}(t)dt
-\frac{1}{4\pi}\int_0^{\infty} H_{21}(t)H_{23}(t)dt
  \end{eqnarray*}
where $\int_0^{\infty} H_{21}(t)dt=.10652624145$ and
$\int_0^{\infty} H_{21}(t)H_{23}(t)dt=.090270992310$. Also
$\int_0^{\infty} H_{21}(t)H_{22}(t)dt=0.3240703564$  and
$\int_0^{\infty} H_{21}(t)H_{23}(t)H_{24}(t)dt=0.2243026030$

  For $u=O(n^{5/4})$ and  $0<x<1$, let
$x=1-\frac{t}{n+t}$. Thus  by (1.6), $EM_u(0,1
)=\left(\frac{n}{(n+t)^2} \right) \int_{0}^{\infty} f_n
\left(1-\frac{t}{n+t}\right)dt,$ where by (1.5) and (1.6) we have
that
$$
\frac{n}{(n+t)^2}G_{1}\bigg(\!1\!-\!
\frac{t}{n+t}\bigg)\!=\!H_{31}(t)\!+\!O(n^{\!-1}),\;\;\;\;
G_{2}\left(\!1\!-\!\frac{t}{n+t}\right)\!=\!\frac{2u}{n^{3/2}\sqrt{\pi}}H_{32}(t)\!
+\!O(n^{\!-5/4}),
$$
$$G_{3}\left(1-\frac{t}{n+t}\right)=H_{33}(t)+O(n^{-1}),\;\;\;\;\;
G_{5}\left(1-\frac{t}{n+t}\right)=1+O(n^{-1/2}),\hspace{.3in}\eqno(2.3)$$
$$
G_{4}\left(1-\frac{t}{n+t}\right)=\frac{2u}{n^{3/2}\sqrt{\pi}}H_{34}(t)+O(n^{-5/4});$$
where
$$
H_{31}(t)={\frac {-1}{31t}} \bigg( (
-1/4t-5/4{t}^{2}+11{t}^{3}+{t}^{5}+{\frac {5}{16}}+11/4{t}^{4} )
{e ^{-4t}}$$$$+ ( -23/2{t}^{3}-1/2{t}^{5}+{\frac {13}{4}}{t}^{4}
+{\frac {65}{4}}{t}^{2}-{\frac {59}{8}}t+{\frac {31}{16}})
{e^{-2t}}$$$$+ ( 12t-4{t}^{3}-18{t}^{2}+3{t}^{4}-9/2
 ) {e^{-3t}}+ ( -{\frac {35}{16}}-3/8t ) {e^{-6
t}}+ $$$$( 5/2{t}^{2}-7/2t+{\frac {17}{4}} ) {e^{-5t}}-1/ 16+(
1/2{t}^{2}-1/2t+1/4 ) {e^{-t}} \bigg)$$$$\times
\bigg(35-(32+32t+160t^2)e^{-t}+(80t^5-216t^4+600t^3+324t^2-588t-294)e^{-2t}$$$$
+(544+1632t+1056t^2-512t^3-192t^4)e^{-3t}-(253+1012t+1400t^2+736t^3+172t^4+16t^5)e^{-4t}
\bigg)^{1/2}$$$$\times \bigg( ( -{\frac {2161}{124}}{t}^{5}-{
\frac {129}{31}}{t}^{7}+{\frac {1949}{124}}{t}^{6}-{\frac {1043}{
124}}{t}^{3}+{t}^{8}+{\frac {3717}{992}}-{\frac {1521}{248}}t-{
\frac {32777}{496}}{t}^{2}+{\frac {501}{16}}{t}^{4}) {e^{-4
t}}$$$$+ ( {\frac {1073}{248}}{t}^{3}+{\frac {5}{31}}{t}^{5}-{
\frac {73}{124}}{t}^{4}-{\frac {523}{496}}{t}^{2}+{\frac {485}{248
}}t-{\frac {733}{496}} ) {e^{-2t}}$$$$+ ( -{\frac {232}{31}
}{t}^{4}+{\frac {375}{248}}t-{\frac {364}{31}}{t}^{3}+{\frac {
486}{31}}\,{t}^{2}+{\frac {1043}{248}}-{\frac {34}{31}}\,{t}^{6}+{
\frac {42}{31}}{t}^{5} ) {e^{-3t}}$$$$+ ( {\frac {21891}{
248}}t+{\frac {20383}{496}}-{\frac {8}{31}}{t}^{8}-{\frac {68}{31}
}{t}^{7}-{\frac {3907}{248}}{t}^{4}-{\frac {547}{62}}{t}^{6}-{
\frac {29851}{248}}{t}^{3}-{\frac {12069}{496}}{t}^{2}-{\frac {
1191}{124}}{t}^{5} ) {e^{-6\,t}}$$$$+ ( {\frac {12731}{124}}
{t}^{2}-{\frac {7153}{248}}t-{\frac {6887}{248}}+{\frac {1697}{62}
}{t}^{5}-{\frac {2335}{62}}{t}^{4}+{\frac {5627}{62}}{t}^{3}-{
\frac {173}{31}}{t}^{6}-{\frac {78}{31}}{t}^{7}) {e^{-5t} }$$$$+(
-{\frac {35}{124}}{t}^{2}-{\frac {111}{248}}t-{\frac {21 }{248}} )
{e^{-t}}$$$$+ ( {\frac {4155}{496}}{t}^{4}+{\frac {
6281}{248}}t+{\frac {3}{31}}{t}^{6}+{\frac {787}{31}}{t}^{3}+{
\frac {12155}{1984}}+{\frac {19097}{496}}{t}^{2}+{\frac
{87}{62}}{ t}^{5} ) {e^{-8t}}$$$$+{\frac {115}{1984}}+ \bigg(
{\frac {1243}{ 62}}{t}^{3}+{\frac {18}{31}}{t}^{6}-{\frac
{20247}{248}}t-{ \frac {2018}{31}}{t}^{2}+{\frac
{401}{62}}{t}^{5}-{\frac {6423}{ 248}}+{\frac {51}{2}}{t}^{4}
\bigg) {e^{-7t}} \bigg) ^{-1},
$$
$$
H_{32}(t)=\sqrt {160}{t}^{3/2} \bigg(( {\frac {253}{40}}t+{\frac {
23}{5}}{t}^{3}+{\frac {35}{4}}{t}^{2}+{\frac {253}{160}}+{\frac {
43}{40}}{t}^{4}+1/10{t}^{5} ) {e^{-4t}}$$$$+ ( {\frac {
147}{80}}+{\frac {27}{20}}{t}^{4}-{\frac {81}{40}}{t}^{2}-1/2{t}
^{5}+{\frac {147}{40}}t-{\frac {15}{4}}\,{t}^{3} ) {e^{-2t}}$$$$
+ ( {\frac {16}{5}}{t}^{3}+6/5{t}^{4}-{\frac {51}{5}}t-{ \frac
{33}{5}}{t}^{2}-{\frac {17}{5}} ) {e^{-3t}}-{\frac {7} {32}}+ (
1/5+1/5t+{t}^{2}) {e^{-t}} \bigg)^ {1/2}
$$$$ \times \bigg(
 ( 16{t}^{5}-4t+176{t}^{3}-20{t}^{2}+44{t}^{4}+5
 ) {e^{-4t}}\hspace{2in}$$$$+ ( 260{t}^{2}+31+52{t}^{4}-184{t}^{3}-
118t-8{t}^{5} ) {e^{-2t}}$$$$+ ( -72-288{t}^{2}+192t+
48{t}^{4}-64{t}^{3} ) {e^{-3t}}+ ( -8t+4+8{t}^{2}
 ) {e^{-t}}$$$$\hspace{2in}+ ( 68-56t+40{t}^{2}t) {e^{-5t}}-1-6
{e^{-6t}}t-35{e^{-6t}} \bigg) ^{-1/2}
$$

and
$$
H_{33}(t)=1/20\,\sqrt {20}\bigg( \bigg(  \left(
13/2{t}^{3}-3/2t-{\frac {75 }{4}}{t}^{2}+{\frac {39}{4}}-3{t}^{4}
\right) {e^{-2t}}\hspace{1in}$$$$+ \left( {t}^{3}+{\frac
{29}{4}}{t}^{2}+{\frac {55}{8}}+{\frac {33}{2}}t
 \right) {e^{-4t}}+ \left( -{\frac {39}{2}}t+15{t}^{2}-{\frac {
31}{2}}+3{t}^{3} \right) {e^{-3t}}-5/8$$$$+ \left( -1/2+9/2t+{t}^{2}
 \right) {e^{-t}} \bigg) ^{2} \bigg(  \left( 9/4-9/2t+3/2{t}^{2}-
{t}^{3} \right)
{e^{-2t}}$$$$+(3t-4)e^{-3t}+(15/8+1/2t)e^{-4t}+te^{-t}-1/8 \bigg)
^{-1}$$$$\times
 \bigg(  \left( {\frac {27}{20}}{t}^{4}-{\frac {15}{4}}{t}^{3}-{
\frac {81}{40}}{t}^{2}+{\frac {147}{80}}-1/2{t}^{5}+{\frac {147}{
40}}t \right) {e^{-2t}}$$$$+ \left( {\frac {23}{5}}{t}^{3}+{\frac {
253}{160}}+{\frac {43}{40}}{t}^{4}+1/10{t}^{5}+{\frac {253}{40}}
t+{\frac {35}{4}}{t}^{2} \right) {e^{-4t}}$$$$+ \left( 6/5{t}^{4}-{
\frac {17}{5}}-{\frac {51}{5}}t-{\frac {33}{5}}{t}^{2}+{\frac {16}
{5}}{t}^{3} \right) {e^{-3t}}-{\frac {7}{32}}+ \left( {t}^{2}+1/5
t+1/5 \right) {e^{-t}} \bigg) ^{-1}\bigg)^{1/2}
$$
$$
H_{34}(t)={\frac {1}{\sqrt{64}}} \bigg( 5+ \left( 4-36t-8{t}^{2}
 \right) {e^{-t}}+ \left( 24{t}^{4}-52{t}^{3}+150{t}^{2}+12t-
78 \right) {e^{-2t}}$$$$+ \left( 156t-24{t}^{3}-120{t}^{2}+124
 \right) {e^{-3t}}- \left( 8{t}^{3}+58{t}^{2}+132t+55 \right)
{e^{-4t}} \bigg) {t}^{3/2}$$$$\times \bigg( \bigg( \left( -{ \frac
{65}{2}}{t}^{2}+{\frac {59}{4}}t+23{t}^{3}-{\frac {31}{8}}
+{t}^{5}-13/2{t}^{4} \right) {e^{-2t}}$$$$+ \left( 1/2t-22{t}^{3}+
5/2{t}^{2}-11/2{t}^{4}-2{t}^{5}-5/8 \right) {e^{-4t}}+ \left(
-24t+8{t}^{3}+9+36{t}^{2}-6{t}^{4} \right) {e^{-3t}}$$$$+
 \left( -{t}^{2}-1/2+t \right) {e^{-t}}+ \left( -5{t}^{2}-17/2+7t
 \right) {e^{-5t}}+1/8+(35/8+3/4t)e^{-6t} \bigg) $$$$\times \bigg(  \left( {t}^{3}+9/2\,t-9/4-3/2\,{t}^{2}
\right) {e^{ -2t}}-te^{-t}+(4-3t)e^{-3t}-(15/8+1/2t)e^{-4t}+1/8
\bigg)\bigg)^{-1/2}
$$

As $t\rightarrow \infty $ we have that
$$H_{31}(t)\sim \frac{4\sqrt{35}}{115t},\;\;\;\;\; H_{33}(t)\sim
\frac{5}{\sqrt{35}},\;\;\;\;\; H_{32}(t)\sim  \sqrt{35}\,
t^{3/2},\;\;\;\;\; H_{34}(t)\sim 5t^{3/2}.
$$
For any real numbers $A$ and $B$ we have that
$$
\frac{A}{t}-\frac{B\sqrt{t}}{n^{3/2}}=\frac{A}{t}
-\frac{B\sqrt{t}}{n^{3/2}+(B/A)t^{3/2}}+O(n^{-3}).\eqno(2.4)
$$
Let   $\; a=\frac{\sqrt{35}}{115\pi}$ and $\;
b=\frac{10u}{23\pi^{3/2}}$, $\; c=\frac{1}{23\pi}$ and $\;
d=\frac{14u}{23\pi^{3/2}}$. Now by (2.3), (2.4), and above
calculations we have that
\begin{eqnarray*}
EM_u(0,1)\!&\!=\!&\!\frac{n}{(n+t)^2}\int_0^{\infty}f_n(1-\frac{t}{n+t})dt=
\frac{1}{4\pi}\int_0^{\infty}
\left(H_{31}(t)-\frac{4\sqrt{35}I_{[t\geq 1]} }{115t}\right)dt\\&&
+\frac{u}{2(n\pi)^{3/2}}\int_0^{\infty}\left(H_{31}(t)H_{32}(t)-\frac{28\sqrt{t}I_{[t\geq
1]}}{23}\right)dt\\&&
-\frac{1}{4\pi}\int_0^{\infty}\left(H_{31}(t)H_{33}(t)-\frac{4I_{[t\geq
1]}}{23t}\right)dt\\&&
-\frac{u}{2(n\pi)^{3/2}}\int_0^{\infty}\left(H_{31}(t)H_{33}(t)H_{34}(t)
-\frac{20\sqrt{t}I_{[t\geq
1]}}{23}\right)dt\\&&+\int_1^{\infty}\frac{a}{t}-\frac{b\sqrt{t}}{n^{3/2}+(b/a)t^{3/2}}dt
-\int_1^{\infty}\frac{c}{t}-\frac{d\sqrt{t}}{n^{3/2}+(d/c)t^{3/2}}dt
+O(n^{-3})
\end{eqnarray*}
where$$\int_0^{\infty}\!\!\left(\!H_{31}(t)\!-\!\frac{4\sqrt{35}I_{[t\geq
1]}}{115t}\! \right)\!dt\!=\!-0.2545810,\;
\int_0^{\infty}\!\!\left(\!H_{11}(t)H_{33}(t)\!-\!\frac{4I_{[t\geq
1]}}{23t}\! \right)\!dt\!=\!-0.2085374$$ Also
 $$\int_0^{\infty}\left(H_{31}(t)H_{32}(t)
 -\frac{28\sqrt{t}I_{[t\geq 1]}}{23}\right)dt=-4.808177963,$$
$$\int_0^{\infty}\left(H_{31}(t)H_{33}(t)H_{34}(t)-\frac{20\sqrt{t}I_{[t\geq
1]}}{23}\right)dt=-2.774789804$$ and by the above assumption for
$\,a$, $\,b,\,$ $\,c,\,$ and $\, d$ we have that
$$\int_1^{\infty}\frac{a}{t}-\frac{b\sqrt{t}}{n^{3/2}+\frac{b}{a}t^{3/2}}dt
=\frac{2a}{3}\mbox{ln} \bigg( \frac{a}{b}n^{3/2}+1
\bigg)=\frac{2\sqrt{35}}{345\pi}\mbox{ln}
\bigg(\frac{\sqrt{35\pi}}{50u}n^{3/2}+1 \bigg),
$$
$$\int_1^{\infty}\left(\frac{c}{t}-\frac{d\sqrt{t}}{n^{3/2}+(d/c)t^{3/2}}\right)dt
=\frac{2c}{3}\mbox{ln} \bigg( \frac{c}{d}n^{3/2}+1
\bigg)=\frac{2}{69\pi}\mbox{ln}
\bigg(\frac{\sqrt{\pi}}{14u}n^{3/2}+1 \bigg)
$$

 For $u=O(n^{1/4})$ and $-1<x<0$, let $x=-1+\frac{t}{n+t}$.
Thus  by (1.6), $EM_u(-1,0 )=\left( \frac{n}{(n+t)^2}\right )
\int_{0}^{\infty} f_n\left(-1+\frac{t}{n+t}\right)dt,\;$ where by
(1.6), (1.5) we have that
$$
\frac{n}{(n+t)^2}G_{1}\left(\!-\!1\!+\!\frac{t}{n+t}\right)\!=\!
H_{41}(t)+O(n^{\!-1}),\;\;\;
G_{2}\left(\!-\!1\!+\!\frac{t}{n+t}\right)\!=\!
\frac{2u}{\sqrt{n\pi}}H_{42}(t)\!+\!O(n^{\!-5/4}),$$
$$\;\;
G_{3}\left(\!-\!1\!+\!\frac{t}{n+t}\right)\!=\!H_{43}(t)\!+\!O(n^{\!-1}),
\;\;\;\;G_{5}\left(\!-1\!+\!\frac{t}{n+t}\right)\!=\!1
+\!O(n^{\!-1/2}),\eqno(2.5)
 $$
$$
 G_{4}\left(-1+\frac{t}{n+t}\right)= \frac{2u}{\sqrt{n\pi}}H_{44}(t)+O(n^{-5/4}),$$
where
$$
H_{41}(t)=\frac{1}{16t}  \bigg(
 ( 2t^5-{\frac {31}{2}}{t}^{4}-{\frac {9}{8}}-{\frac {
9}{2}}t+\frac{3}{2}{t}^{2}+12{t}^{3} ) {e^{-4t}}+ ( 7{t
}^{3}-{\frac {9}{2}}{t}^{2}+{t}^{5}-{\frac {9}{2}}{t}^{4}+{ \frac
{3}{8}}+\frac{3}{4}t ) {e^{-2t}}$$$$+\frac{1}{8}+ ( 6{t}^{2}-
4{t}^{5}-8{t}^{3}-11{t}^{4}+{\frac {15}{4}}t+{\frac {5}{8}}
 ) {e^{-6t}} \bigg)$$$$ \times \bigg(\bigg(3+(12t-32t^3+24t^2-44t^4-16t^5+3)e^{-4t}
 +(16t^5-56t^4+56t^3-12t^2-12t-6)e^{-2t}\bigg)\bigg)^{1/2}\;\; $$$$\times\bigg(
( {\frac {47}{32}} {t}^{3}+{\frac {15}{128}}{t}^{2}-{\frac
{205}{32}}{t}^{5} +{ \frac {257}{32}}{t}^{6}-{\frac
{39}{256}}-{\frac {35}{8}}{t}^{7}+{ \frac {7}{8}}{t}^{8}-{\frac
{39}{64}}t+{\frac {35}{128}}{t}^{4}
 ) {e^{-4t}}$$$$+ ( {\frac {45}{64}}{t}^{3}-{\frac {55}{128
}}{t}^{2}+1/8{t}^{5}-{\frac {15}{32}}{t}^{4}+{\frac {1}{64}}t+
{\frac {1}{128}} ) {e^{-2t}}$$$$+ ( {\frac {279}{128}}{t}^{
2}+{\frac {73}{32}}{t}^{5}+{\frac {25}{128}}-{\frac {79}{64}}{t}^{
3}-{\frac {459}{64}}{t}^{4}+{\frac {75}{64}}t+{t}^{8}+{\frac {165}
{16}}{t}^{6}-9{t}^{7} ) {e^{-6t}}$$$$+{\frac {11}{512}}+
 ( {\frac {507}{128}}{t}^{4}+1/8{t}^{6}+{\frac {147}{16}}{t
}^{5}-{\frac {239}{128}}{t}^{2}-{\frac {37}{512}}-9/2{t}^{7}-{ \frac
{37}{64}}t-{\frac {27}{16}}{t}^{3}-2{t}^{8} ) {e^{-8 t}} \bigg)
^{-1} ,
$$
and
$$
H_{42}(t)=2\sqrt{t}\,\bigg(
(6+12t+12t^2-56t^3+56t^4-16t^5)e^{-2t}$$$$+(16t^5+44t^4+32t^3-24t^2-12t-3)e^{-4t}-3\bigg)^{1/2}$$$$\times
\bigg[(36t^4-8t^5-56t^3+36t^2-6t-3)e^{-2t}+(9+36t-12t^2-96t^3+124t^4-16t^5)e^{-4t}$$
$$+(32t^5+88t^4+64t^3-48t^2-30t-5)e^{-6t}-1\bigg]^{-1/2}
,
$$
$$
H_{43}(t)=8\bigg( ( -1/4-1/2t-5/2{t}^{3}+{t}^{4}+7/4{t}^{2}
 )e^{-2t}+1/8+( 1/2t-5/4t^{2}-{t}^{3}+1/8
 )e^{-4t}\bigg)$$$$\times \bigg[\bigg(3+(16t^5-56t^4+56t^3-12t^2-12t-6)e^{-2t}
 +(3+12t+24t^2-32t^3-44t^4-16t^5)e^{-4t}\bigg)$$$$\times
\bigg(1+(2+4t-12t^2+8t^3)e^{-2t}+(16t^3-8t^2-12t-3)e^{-4t}\bigg)\bigg]^{-1/2}
$$
$$
H_{44}(t)=2\sqrt{t}\,\bigg(
1+(8t^4-20t^3+14t^2-4t-2)e^{-2t}+(4t-10t^2-8t^3+1)e^{-4t}
\bigg)$$$$\times\bigg[\bigg(
(3+56t^3+8t^5-36t^4+6t-36t^2)e^{-2t}+(16t^5+12t^2-9+96t^3-36t-124t^4)e^{-4t}$$$$+(30t-88t^4+5+48t^2-32t^5-64t^3)e^{-6t}+1
\bigg)
$$$$\times \bigg( (3+12t+8t^2-16t^3)e^{-4t}-(2+4t-12t^2+8t^3)e^{-2t}-1 \bigg) \bigg]^{-1/2}$$

As $t\rightarrow \infty $ we have that
$$H_{41}(t)\sim \frac{4\sqrt{3}}{11t},\;\;\;\;\; H_{43}(t)\sim
\frac{1}{\sqrt{3}},\;\;\;\;\; H_{42}(t)\sim 2\sqrt{3t} ,\;\;\;\;\;
H_{44}(t)\sim 2\sqrt{t}.
$$
 For any real numbers $A$ and $B$ we have that
$$
\frac{A}{t}-\frac{B/\sqrt{t}}{n^{1/2}}=\frac{A}{t}
-\frac{B/\sqrt{t}}{n^{1/2}+(B/A)t^{1/2}}+O(n^{-1}). \eqno(2.6)
$$
Let   $\; a=\frac{\sqrt{3}}{11\pi}$ and $\;
b=\frac{4u}{11\pi^{3/2}}$, $\; c=\frac{1}{11\pi}$ and $\;
d=\frac{12u}{11\pi^{3/2}}$.  Now by (2.5), (2.6) and above
calculations we have that
\begin{eqnarray*}
EM_u(-1,0)=\!&\!\!&\!\hspace{-.2in}\frac{n}{(n+t)^2}\int_0^{\infty}f_n(-1+\frac{t}{n+t})dt=
\frac{1}{4\pi}\int_0^{\infty}\left(\!
H_{41}(t)-\frac{4\sqrt{3}I_{[t\geq
1]}}{11t}\right)dt\\&&\hspace{-.5in}+\!
\frac{u}{2\pi\sqrt{n\pi}}\int_0^{\infty}\!\left(\!
H_{41}(t)H_{42}(t)\!-\!\frac{24I_{[t\geq
1]}}{11\sqrt{t}}\right)dt\!-\!
\frac{1}{4\pi}\int_0^{\infty}\!\left(
H_{41}(t)H_{43}(t)\!-\!\frac{4I_{[t\geq
1]}}{11t}\right)dt\\&&\!\!\!\!\!\!\hspace{-.1in}-
\frac{u}{2\pi\sqrt{n\pi}}\int_0^{\infty}\left(
H_{41}(t)H_{43}(t)H_{44}(t)-\frac{8I_{[t\geq
1]}}{11\sqrt{t}}\right)dt
\\&&\!\!\!\!\!\!\hspace{-.1in}+\int_1^{\infty}\frac{a}{t}-\frac{b/\sqrt{t}}{n^{1/2}+(b/a)t^{1/2}}dt
-\int_1^{\infty}\frac{c}{t}-\frac{d/\sqrt{t}}{n^{1/2}+(d/c)t^{1/2}}dt
+O(n^{-1})
\end{eqnarray*}
where  $\int_0^{\infty}\left( H_{41}(t)-\frac{4\sqrt{3}I_{[t\geq
1]}}{11t}\right)dt=-0.1146419848,$ and
  $ \int_0^{\infty}\left(
H_{41}(t)H_{43}(t)-\frac{4I_{[t\geq
1]}}{11t}\right)dt=-0.0801100983$. Also

\noindent
 $\int_0^{\infty}\!
\left(\! H_{41}(t)H_{42}(t)\!-\!\frac{24I_{[t\geq
1]}}{11\sqrt{t}}\! \right)dt\!=\!-0.7769335,\;
   \int_0^{\infty}\!\left(\!
H_{41}(t)H_{43}(t)H_{44}(t)\!-\!\frac{8I_{[t\geq
1]}}{11\sqrt{t}}\!\right)dt\!=\!-0.1820104$, and by the above
assumption for $\,a$, $\,b,\,$ $\,c,\,$ and $\, d$ we have that
$$\int_1^{\infty}\frac{a}{t}-\frac{b/\sqrt{t}}{n^{1/2}+\frac{b}{a}t^{1/2}}dt
=2a\mbox{ln} \bigg( \frac{a}{b}n^{1/2}+1
\bigg)=\frac{2\sqrt{3}}{11\pi}\mbox{ln}
\bigg(\frac{\sqrt{3\pi}}{4u}n^{1/2}+1 \bigg),
$$$$\int_1^{\infty}\left(\frac{c}{t}-\frac{d/\sqrt{t}}{n^{1/2}+(d/c)t^{1/2}}\right)dt
=2c\mbox{ln} \bigg( \frac{c}{d}n^{1/2}+1
\bigg)=\frac{2}{11\pi}\mbox{ln}
\bigg(\frac{\sqrt{\pi}}{12u}n^{1/2}+1 \bigg)
$$
Simplifying these calculations lead to the result of the theorem.

\end{document}